\documentclass[12pt,a4paper]{amsart}
\usepackage{amsfonts}
\usepackage{amsthm}
\usepackage{amsmath}
\usepackage{amscd}
\usepackage[latin2]{inputenc}
\usepackage{t1enc}
\usepackage[mathscr]{eucal}
\usepackage{indentfirst}
\usepackage{graphicx}
\usepackage{graphics}
\usepackage{pict2e}
\usepackage{epic}
\usepackage{amssymb}
\usepackage{amstext}
\usepackage{ucs}
\usepackage[plainpages=false,pdfpagelabels,pagebackref]{hyperref}
\usepackage[capitalise]{cleveref}

\numberwithin{equation}{section}
\usepackage[margin=2.9cm]{geometry}
\usepackage{epstopdf}

\theoremstyle{plain}
\newtheorem{Th}{Theorem}[section]
\newtheorem{Lemma}[Th]{Lemma}

 \theoremstyle{definition}

\newtheorem{?}[Th]{Problem}

\begin{document}

\title{ESSENTIAL SPHERICAL ISOMETRIES}

\author{Marcel Scherer}

\begin{abstract} A result due to  Williams, Stampfli and Fillmore shows that an essential isometry $T$ on a Hilbert space $\mathcal{H}$ is a compact perturbation of an isometry if and only if ind$(T)\le 0$. A recent result of S. Chavan yields an analogous characterization of essential spherical isometries $T=(T_1,\dots,T_n)\in\mathcal{B}(\mathcal{H})^n$ with dim($\bigcap_{i=1}^n\ker(T_i))\le$ dim$(\bigcap_{i=1}^n\ker(T_i^*))$. In the present note we show that in dimension $n>1$ the result of  Chavan holds without any condition on the dimensions of the joint kernels of $T$ and $T^*$.
\end{abstract}

\maketitle

\section{Introduction} Let $\mathcal{H}$ be an infinite-dimensional Hilbert space, $n\in\mathbb{N}$ and let $T=(T_1,\dots,T_n)\in\mathcal{B}(\mathcal{H})^n$ be a tuple of continuous linear operators. Depending on the context, we use the notation $T$ both for the tuple $T = (T_1,\dots,T_n)$ and for the column operator $T: H\to H^n,
x \mapsto (T_1(x),\dots,T_n(x))$. The adjoint of the column operator $T$ acts as the
row operator $T^*:\mathcal{H}^n\to\mathcal{H},(x_i)_{i=1}^n\mapsto\sum_{i=1}^n T_i(x_i)$. We say that $T$ is \textit{essentially commuting} if $T_iT_j-T_jT_i$ is compact for all $i,j\in\{1,\dots,n\}$ and we call $T$ an \textit{essential spherical isometry} if $T$ is essentially commuting and its column operator is an essential isometry. For $n=1$, we call $T$ an \textit{essential isometry}.\\
The following theorem of Fillmore, Stampfli and Williams \cite{FSW72} gives the motivation for this paper: 
\begin{Th} Let $T\in\mathcal{B}(\mathcal{H}).$ Then $T$ is an essential isometry with 
$\dim(\ker(T))\le \dim(\ker(T^*))$ if and only if $T$ is a compact perturbation of an isometry.
\label{thm1}
\end{Th}
The \textit{joint kernel} $\ker(T)$ of an operator tuple $T=(T_1,\dots,T_n)$ is defined as $\ker(T)=\bigcap_{i=1}^n\ker(T_i)$. Our aim is to improve the following theorem of Chavan \cite{article}.
\begin{Th} Let $T\in\mathcal{B}(\mathcal{H})^n$ and $\dim(\ker(T))\le \dim(\ker(T^*))$. Then $T$ is a compact perturbation of a tuple $V\in\mathcal{B}(\mathcal{H})^n$ that is essentially commuting and whose associated column operator is an isometry if and only if $T$ is an essential spherical isometry.
\label{thm2}
\end{Th}
For $n=1$ this theorem is exactly \cref{thm1} and it is clear that the condition $\dim(\ker(T))\le \dim(\ker(T^*))$ is necessary because there are essential isometries with $\dim(\ker(T))>\dim(\ker(T^*))$. For example, the left shift $L$ on $\ell^2$ fulfills $\dim(\ker(L))=1>0=\dim(\ker(L^*))$. But is there any essential spherical isometry $T$ with $\dim(\ker(T))>\dim(\bigcap_{i=1}^n\ker(T_i^*))$? The next section will clarify this question.

\section{main result}
First, we recall some properties of the column operator induced by an essential spherical isometry. For a proof, see \cite[Lemma 2.6.13]{book}.
\begin{Lemma}Let $T\in\mathcal{B}(\mathcal{H},\mathcal{H}^n)$. The following conditions are equivalent:
\renewcommand{\labelenumi}{(\roman{enumi})}
\begin{enumerate}
\item There is an operator $S\in\mathcal{B}(\mathcal{H}^n,\mathcal{H})$ such that $id-ST$ is compact.
\item The range of $T$ is closed and the null-space $\ker(T)$ is finite-dimensional.
\end{enumerate}
\label{lem2}
\end{Lemma}
\ \\

If $T$ is an essential spherical isometry, then the column operator of $T$ satisfies the first condition in \cref{lem2}. Hence dim$(\ker(T))<\infty$. Using this observation we can prove the next theorem.
\begin{Th} Let $T=(T_1,\dots,T_n)\in\mathcal{B}(\mathcal{H})^n$ be an essential spherical isometry. If $n>1$, then the kernel of the row operator $T^*:\mathcal{H}^n\to\mathcal{H}$ is infinite-dimensional.
\label{thm3}
\end{Th}
\textit{Proof}.  Assume that $n>1$ and that the kernel of $T^*:\mathcal{H}^n\to\mathcal{H}$ is finite dimensional. Since by \cref{lem2} the column operator $T:\mathcal{H}\to\mathcal{H}^n$ has closed range and finite-dimensional kernel, it follows that $T:\mathcal{H}\to\mathcal{H}^n$  is a Fredholm operator. By Atkinson's theorem there is a row operator $S=(S_1,\dots,S_n):\mathcal{H}^n\to\mathcal{H}$ such that $1_{\mathcal{H}^n}-TS\in\mathcal{B}(\mathcal{H}^n)$ is compact. Let us denote by $\pi_i:\mathcal{H}^n\to\mathcal{H}, (x_1,\dots,x_n)\mapsto x_i$, the canonical projections. For $i,j\in\{1,\dots,n\}$ with $i\neq j$, the operators
\begin{equation*}
-T_iS_j=\pi_i(1_{\mathcal{H}^n}-TS)\pi_j^*,
\end{equation*}
\begin{equation*}
1_{\mathcal{H}}-T_jS_j=\pi_j(1_{\mathcal{H}^n}-TS)\pi_j^*
\end{equation*}
and $T_iT_j-T_jT_i$ are compact. Hence 
\begin{equation*}
T_i=T_i(1_\mathcal{H}-T_jS_j)+(T_iT_j-T_jT_i)S_j+T_jT_iS_j
\end{equation*}
is compact for $i=1,\dots,n$. But then the column operator $T:\mathcal{H}\to\mathcal{H}^n$ is a compact Fredholm operator and hence $\mathcal{H}$ would have to be finite dimensional. This contradiction completes the proof. \qed\\
\ \\

\begin{Th}Let $n>1$ and $T\in\mathcal{B}(\mathcal{H})^n$. Then $T$ is an essential spherical isometry if and only if $T=V+K$ with a tuple $K\in\mathcal{B}(\mathcal{H})^n$ of compact operators and an essentially commuting tuple $V\in\mathcal{B}(\mathcal{H})^n$ such that $V^*V=1_\mathcal{H}$.
\label{thm4}
\end{Th}
\textit{Proof}. Obviously each tuple $T=V+K$ with $V$ and $K$ as in \cref{thm4} is an essential spherical isometry. Conversely, let $T\in\mathcal{B}(\mathcal{H})^n$ be an essential spherical isometry. By \cref{lem2} the column operator $T:\mathcal{H}\to\mathcal{H}^n$ has finite-dimensional kernel and closed range. By \cref{thm3} the kernel of its adjoint $T^*:\mathcal{H}^n\to\mathcal{H}$ is infinite-dimensional. Let $K\in\mathcal{B}(\mathcal{H},\mathcal{H}^n)$ be a bound linear operator with $K=0$ on $\ker(T)^\perp$ such that $K$ induces an isometry $K:\ker(T)\to\ker(T^*)$. Then $\tilde{T}=T+K\in\mathcal{B}(\mathcal{H},\mathcal{H}^n)$ is a compact perturbation of $T\in\mathcal{B}(\mathcal{H},\mathcal{H}^n)$ such that
\begin{equation*}
\lVert \tilde{T}(x+y) \rVert^2=\lVert T(y)+K(x) \rVert^2=\lVert T(y) \rVert^2+\lVert K(x) \rVert^2
\end{equation*}
for $x\in\ker(T)$ and $y\in\ker(T)^\perp$. But then $\tilde{T}$ is bounded below and $A=\tilde{T}^*\tilde{T}\in\mathcal{B}(\mathcal{H})$ is positive and invertible. Let us denote by $\mathcal{B}(\mathcal{H})\to\mathcal{C}(\mathcal{H})=\mathcal{B}(\mathcal{H})/\mathcal{K}(\mathcal{H}), S\mapsto [S]$, the quotient map into the Calkin-algebra of $\mathcal{H}$. By construction $A$ is a compact perturbation of the identity operator. Hence $[A^{-\frac{1}{2}}]^2=[A^{-1}]=[1_\mathcal{H}]$ and, since positive square roots in the $C^*$-algebra $\mathcal{C}(\mathcal{H})$ are uniquely determined, the operator $A^{-\frac{1}{2}}$ is a compact perturbation of the identity operator. Then $V=\tilde{T}A^{-\frac{1}{2}}=(T+K)A^{-\frac{1}{2}}\in\mathcal{B}(\mathcal{H},\mathcal{H}^n)$ is a compact perturbation of $T$ with $V^*V=A^{-\frac{1}{2}}AA^{-\frac{1}{2}}=1_\mathcal{H}$. Thus the proof of \cref{thm4} is complete.
\qed \\
\ \\
Another question formulated in \cite{article} whether if there exists an essential spherical isometry which is no compact perturbation of a commuting tuple. To answer this question we recall the following result from \cite{article3}.
\begin{Th}
Let $A,B\in\mathcal{B}(\mathcal{H})$ with $BA-1_{\mathcal{H}}\in\mathcal{K}(\mathcal{H})$ and $AB-1_{\mathcal{H}}\in\mathcal{K}(\mathcal{H})$. Then there are compact perturbations $\tilde{A}$ of $A$ and $\tilde{B}$ of $B$ with $\tilde{A}\tilde{B}=\tilde{B}\tilde{A}$ if and only if $\mathrm{ind}(A)=\mathrm{ind}(B)=0$.
\label{prop}
\end{Th}

We now choose $\mathcal{H}=\ell^2$ and denote by $R$ be the right shift operator. Then $T=(T_1,T_2)=\frac{1}{\sqrt{2}}(R,R^*)$ is an essential spherical isometry. Let us assume that there are compact perturbations $\tilde{T}_1,\tilde{T}_2$ of $T_1,T_2$ with $\tilde{T}_1\tilde{T}_2=\tilde{T}_2\tilde{T}_1$. Because $R$ and $R^*$ are inverse to each other modulo $\mathcal{K}(\ell^2)$, the operators $\sqrt{2}\tilde{T}_1$ and $\sqrt{2}\tilde{T}_2$ are also inverse modulo $\mathcal{K}(\ell^2)$. Thus \cref{prop} says that ind$(\sqrt{2}T_1)=0$. However $\sqrt{2}T_1$ is a compact perturbation of $\sqrt{2}R$  and therefore index$(\sqrt{2}T_1)=-1$. This example shows that there are essential spherical isometries which are no compact perturbation of a commuting operator tuple.

\end{document}